\documentclass{aims}
\usepackage{amsmath,amssymb,array,bm,graphicx}
  \usepackage{paralist}
  \usepackage{comment}
  
  \usepackage[normalem]{ulem}
\def\by	{{\mathbf y}}
\def\bx	{{\mathbf x}}
\def\bX	{{\mathbf X}}
\def \btheta {\mathbf{\theta}}
\def\bt {{\mathbf t}}
\def\bq {{\mathbf q}}
\def\bR {{\mathbf R}} 
\def\bP {{\mathbf P}}
\def\bH {{\mathbf H}}
\def\gbR {{\hat{\mathbf R}}} 
\def\br {{\mathbf r}}
\def\R  {{\mathbb R}}
\def\bY  {{\mathbf Y}}
\def\MVN {\mathsf{MVN}}
\DeclareMathOperator{\ceil}{ceil}

\newcommand{\E}{\mathsf{E}}

\usepackage{defns} 
\usepackage{mycleveref} 

\usepackage{caption}
\usepackage{subcaption} 
\usepackage{enumitem}

\usepackage[numbers]{natbib}

\def\currentvolume{X}
 \def\currentissue{X}
  \def\currentyear{200X}
   \def\currentmonth{XX}
    \def\ppages{X--XX}

\title[A surrogate-based approach to non-Gaussian data assimilation]
      {A surrogate-based approach to nonlinear, non-Gaussian joint state-parameter data assimilation}

\author[John Maclean and Elaine T. Spiller]{}


 \email{john.maclean@adelaide.edu.au}
 \email{elaine.spiller@marquette.edu}

\thanks{The first author is supported by the ARC grant DP180100050, and acknowledges past support from ONR grant N00014-18-1-2204. The second author is supported by NSF grant DMS-1821338}

\begin{document}
\maketitle

\centerline{\scshape John Maclean$^*$}
\medskip
{\footnotesize
 \centerline{School of Mathematical Sciences}
   \centerline{University of Adelaide}
   \centerline{SA 5005, Australia}
} 

\medskip

\centerline{\scshape Elaine T. Spiller}
\medskip
{\footnotesize
 \centerline{Department of Mathematical and Statistical Sciences}
   \centerline{Marquette University}
   \centerline{Milwaukee, WI 53201, USA}
}

\bigskip

 \centerline{(Communicated by the associate editor name)}

\begin{abstract}
Many recent advances in sequential assimilation of data into nonlinear high-dimensional models are modifications to particle filters which employ efficient searches of a high-dimensional state space. In this work, we present a complementary strategy that combines 
statistical emulators and particle filters. The emulators are used to learn and offer a computationally cheap approximation to the forward dynamic mapping. This emulator-particle filter (Emu-PF) approach requires a modest number of forward-model runs, but yields well-resolved posterior distributions even in non-Gaussian cases. We explore several modifications to the Emu-PF that utilize mechanisms for dimension reduction to efficiently fit the statistical emulator, and present a series of simulation experiments on an atypical Lorenz-96 system to demonstrate their performance. We conclude with a discussion on how the Emu-PF can be paired with modern particle filtering algorithms.
\end{abstract}

\section{Introduction/Motivation }
 
Data assimilation (DA) -- the process of updating models with data to give state estimates and forecasts complete with attendant uncertainties -- has progressed tremendously over the last three decades \cite{book:evensen,law2015data,reich2015probabilistic,carrassi2018data}.  Sequential data assimilation techniques, those that update current state estimates and forecasts on the fly as data becomes available, fall into two general categories:  Kalman-type filters (KF) and particle-type filters (PF). There are two dominant challenges in sequential DA, namely systems with high-dimensional state spaces and strong non linearity/non Gaussianity. Typically ensembled-based  KF techniques, that use relatively few model runs, can address the former while particle-based techniques that require many model runs can address the latter.

Over roughly the same time frame, the field of statistical surrogates of complex computer models emerged \cite{sacks1989design,Welc:Buck:etal:1992}. Statistical surrogates offer a rapid approximation to the input/output mapping of a computer model based on a typically modest number of training runs. Further, statistical surrogates interpolate available model runs and offer an approximation with built-in uncertainty estimates for utilizing the approximation. In this work we develop particle filters that employ statistical surrogates to approximate mappings of system dynamics.

Often dynamic model forecasts are the computational bottleneck in sequential data assimilation.
Our approach employs Gaussian process emulators (GPs) to learn the mapping from state and/or parameter values at one observation instance to the next. This mapping provides an effective interpolation between model forecasts and makes available slews of additional approximate model forecasts with negligible additional computational burden.  Thus GP emulators are natural to pair with ``sample hungry" DA techniques like particle filters. As such, we can produce finely-sampled non-Gaussian posterior estimates with a modest number of model runs typical of ensemble KF techniques.

Several papers amount to a recent flurry of activity combining modeling learning and data assimilation, each approach with advantages and drawbacks.  \citep{Cleary2020calibrate} combines statistical emulation and data assimilation to aid in model calibration, effectively a smoothing problem. This approach is quite similar to \cite{bayarri2007computer}, but cleverly utilizes an ensemble Kalman filter (EnKF) sequentially between observations to choose a good design (e.g. a well-chosen set of training runs) to fit the statistical emulator.  That emulator then replaces forward model evaluations in MCMC evaluations in the calibration problem. \citep{Bocquet2020bayesian} 
sets up methodology for combining Neural networks (NN) to learn an approximation to the dynamic forward model from noisy observations in a data assimilation framework. Effectively they construct a posterior distribution as one would in a sequential  DA  problem and use expectation maximization algorithm (EM) to compute mode posterior estimates of the NN parameters. This approach is quite appealing as it does not rely on physics-based model for forward propagation, yet simulation experiments in \citep{Brajard2020combining} indicate the need for a significant amount of data to train the NNs sufficiently. \cite{gott:reic20arxiv} is similar in spirit, but uses a random feature map instead of NNs and  they combine their ``physics-agnostic" forward model with an EnKF. Like \cite{Bocquet2020bayesian,Brajard2020combining}, this methodology requires a significant amount of training data.

\begin{figure}
      \centering
      \begin{subfigure}[b]{0.3\textwidth}
          \centering
          \includegraphics[width=\textwidth]{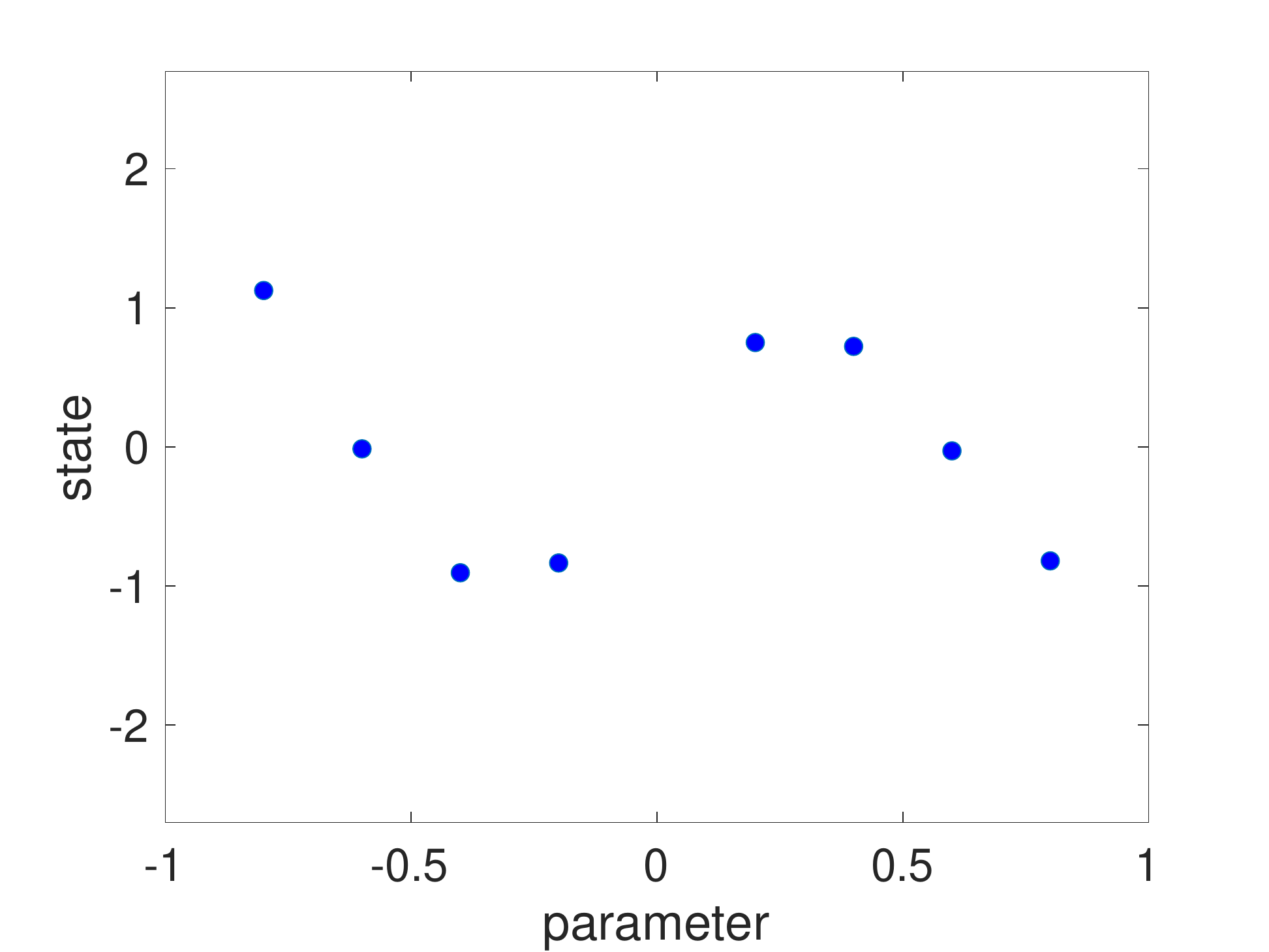}
          \caption{}
          \label{ex1d}
      \end{subfigure}
      \hfill
      \begin{subfigure}[b]{0.3\textwidth}
          \centering
          \includegraphics[width=\textwidth]{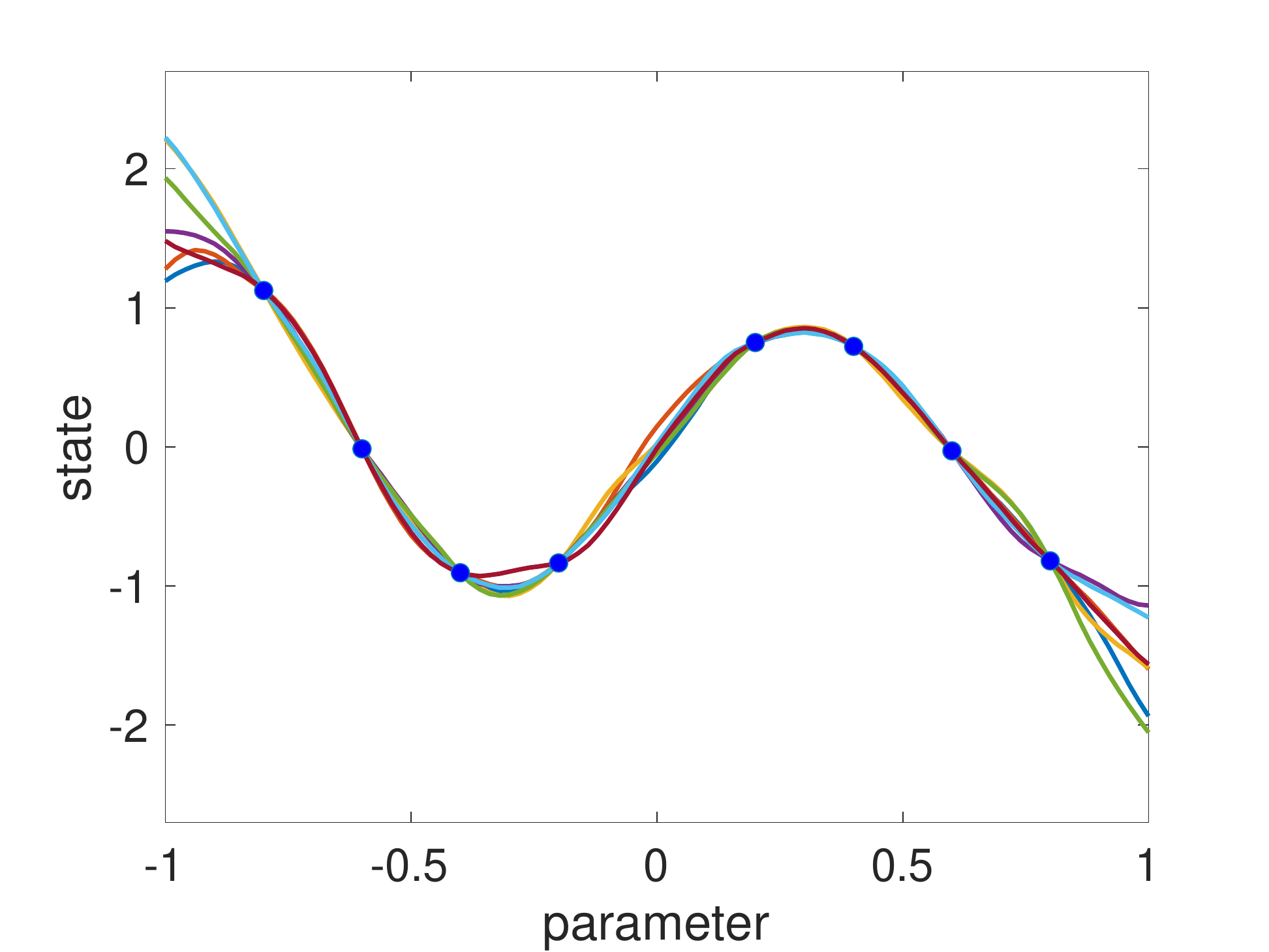}
          \caption{}
          \label{exInterp}
      \end{subfigure}
      \hfill
      \begin{subfigure}[b]{0.3\textwidth}
          \centering
          \includegraphics[width=\textwidth]{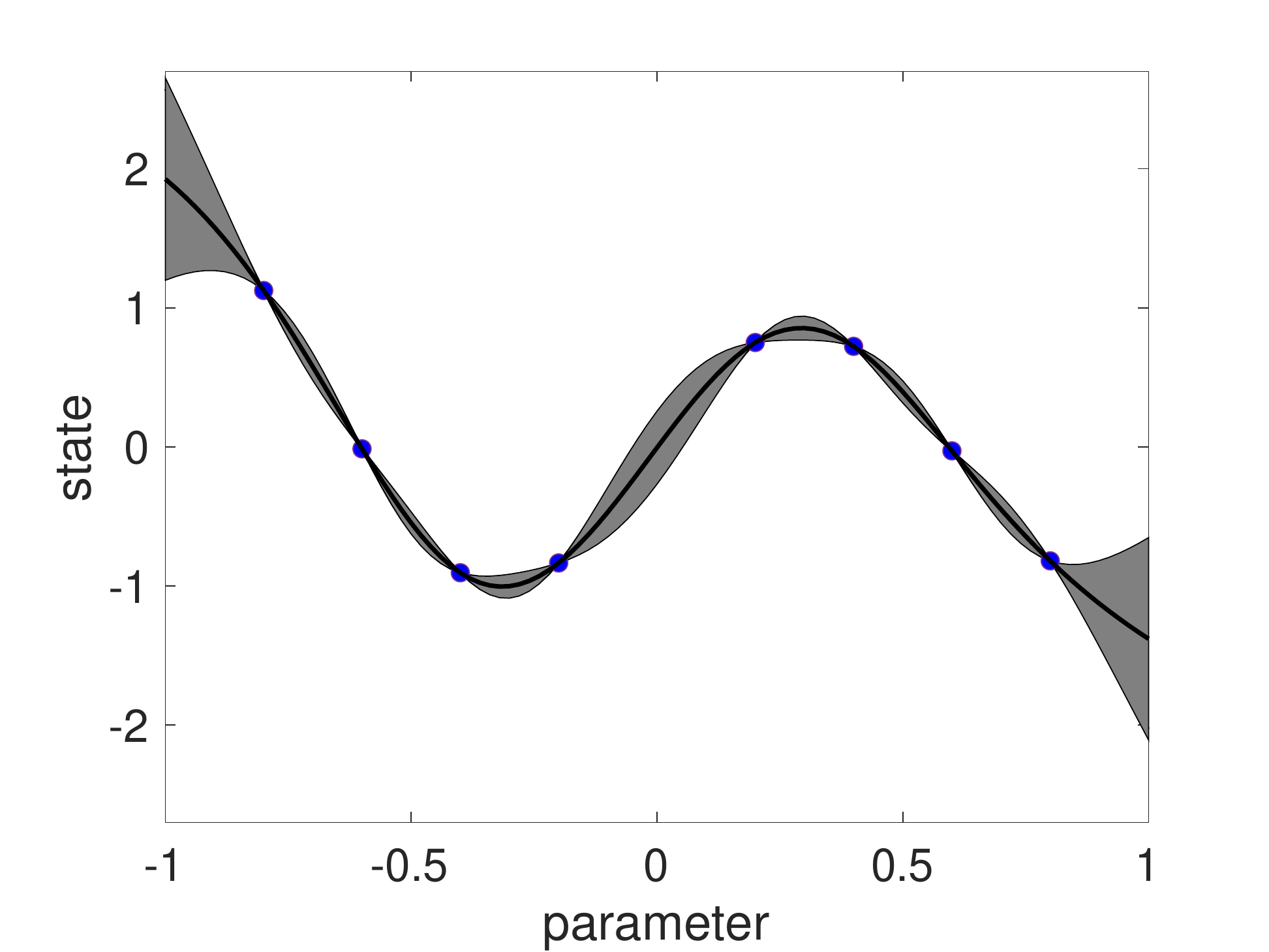}
          \caption{}
          \label{exModel}
      \end{subfigure}
         \caption{Schematic for state dependence on parameters: we plot the state at eight different samples \eqref{ex1d}, then apply a variety of interpolating schemes \eqref{exInterp} and lastly a statistical surrogate \eqref{exModel}. The shaded region in the rightmost plot shows one standard deviation in uncertainty. The second and third plot allow for the state to be estimated at a variety of parameter values.}
         \label{fig:coarseInterp}
 \end{figure}

The core idea of this paper is that interpolation between model forecasts, thought of as functions of the parameter and/or previous state values at a fixed time, may be used to produce additional forecasts, and thus provide a cheap means to improve PF performance. In  the former case, the interpolation exploits smoothness of the state with respect to parameter values, that is not used (nor required) in the usual formulations of the PF. Consider the following pedagogical example in terms of parameter dependence (dependence on previous state values or both is similar in spirit.) Suppose we have some computational model that provides forecasts of a single state variable, and which depends on a single parameter. Notionally this model is expensive to run, and we sample the model forecast at eight parameter values, as in Figure~\ref{ex1d}. Provided the sampling is space filling in parameter space, one might be able to predict the state values at other choices of the parameters, as in Figure~\ref{exInterp}. Using a statistical surrogate instead of a deterministic interpolant, we can further capture the \emph{uncertainty} in the state-parameter dependence at parameter values that have not been sampled, as in Figure~\ref{exModel}. 
These interpolating schemes allow for the output of the computational model to be estimated at many more sampling points than the initial eight. If these fine samples are used as ``model forecasts'' in a PF, then the effect is to obtain a dense estimate of the posterior using relatively few model forecasts.


 
The performance of this method relies on an efficient implementation of the simple interpolate-and-sample concept above. In particular the parameter values at which the model is evaluated are not all fixed, but updated at each observation time, and the interpolating method should be a statistical surrogate that captures uncertainty. These foundational concepts from Data Assimilation and Uncertainty Quantification are introduced in Section~\ref{sec:background}. The surrogate DA scheme is described in Section~\ref{sec:method}, and numerous visualisations of the internal mechanisms and error statistics of the new scheme are contained in Section~\ref{sec:numerics}.

\section{Background}
\label{sec:background}
\subsection{Sequential data assimilation}\label{sec:da}
Let us begin by reviewing the setup for sequential data assimilation and two ``standard" techniques: particle filters (PF) and ensemble Kalman filters (EnKF). In these approaches data or observations of a system are assimilated into a model describing the dynamics of the underlying system to offer an estimate of the state, and of parameters of interest along with attendant uncertainties. This is often achieved by employing Bayes theorem, 
$p(\bx,\btheta \mid \by)=p(\by \mid \bx, \btheta)p(\bx, \btheta)/p(\by)$ where $\bx \in \mathbb{R}^n$ is the state variable, $\by \in \mathbb{R}^m$ is the data, and $\btheta \in \mathbb{R}^p$ are parameters. 

In the sequential case where data are available as a time series, we will follow the notation of Doucet {\em et.~al.}~\citep{SMCM_Doucet}. For observations available at times $\bt=\{t_j,t_{j+1},\ldots,t_k\}$ we use the shorthand $\by_{j:k}=\{\by_j,\ldots,\by_k\}$ and likewise for state variables and parameters at times $\bt$, $\bx_{j:k}=\{\bx_j,\ldots, \bx_k\}$, $\theta_{j:k}=\{\theta_j,\ldots, \theta_k\}$. Ultimately for sequential state-parameter data assimilation, we are interested in describing the posterior distribution 
\begin{align*}
    p(\bx_{0:k}, \btheta_{0:k} \mid \by_{1:k})=\frac{p(\by_{1:k}\mid \bx_{0:k}, \btheta_{0:k})p(\bx_{0:k},\btheta_{0:k})}{p(\by_{1:k})},
\end{align*}
where the marginal distribution in the denominator is given by~$p(\by_{1:k})=\E_{p(\bx_{0:k},\btheta)}[p(\by_{1:k}\mid \bx_{0:k}, \btheta_{0:k})]$. All of these distributions are updated sequentially as data become available. For both PFs and EnKFs, we can think of the representation of the prior and posterior probability density functions as a collection of $N$ ``particles", e.g. state variable and parameter estimates with weights, $\{\bx_j^i, \ \btheta_j^i, \ w_j^i\}_{i=1}^N$. In both cases, the particle states are advanced via system dynamics, e.g. $\dot \bx = \mathcal{M}(\bx,\btheta)$,  from time $t_j$ to time $t_{j+1}$, according to some map
\begin{align} \label{model}
\bx_{j}^i = \varphi(\bx_{j-1}^i,\theta_{j-1}^i) =\bx_{j-1}^i +\int_{t_{j-1}}^{t_{j}} \mathcal{M} (\bx^i, \btheta^i) dt\; . 
\end{align}
This step is typically the most computationally intensive part of the Particle Filter, particularly in the diverse applications in which the dynamical system $\mathcal{M}$ is high-dimensional.

\subsubsection*{Particle filters} For generic particle filters, the particle representations of probability density functions (pdfs) are updated by adjusting the weights via the likelihood as current data are incorporated while state values remain unchanged. Often sequential importance resampling (SIR) \citep{JASA93_443p1032} is employed to overcome inherent filter degeneracy. The idea behind SIR particle filters is to use the posterior distribution from one time step as the prior distribution for the next (along with the state updated by the forward dynamics in \cref{model}) as
\begin{align} \label{bayes}
 p(\bx_{j}, \btheta_{j} \mid \bx_{j-1}, \btheta_{j-1}, \by_{1:j})=\frac{p(\by_{j}\mid \bx_{j})p(\bx_{j},\btheta_j \mid \bx_{j-1}, \btheta_{j-1}, \by_{1:j-1})}{p(\by_{j})}.
\end{align}

Implementing this approach includes updating the weights, with the $j^{th}$ observation yielding
\begin{align} \label{w}
    w_j^i=\frac{p(\by_j \mid \bx^i_j)}{\sum_{i=1}^N p(\by_j \mid \bx^i_j)} w_{j-1}^i
\end{align}
with the posterior approximated by
\begin{align*}
   p(\bx_{j}, \btheta_j \mid \by_{1:j}) \approx \sum_{i=1}^N w_j^i \delta(\bx-\bx_j^i) \delta(\btheta -\btheta_j^i).
\end{align*}
Generally for particle filters, resampling will need to be employed when the effective number of particles, $N_{eff} \approx \sum_{i=1}^N 1/{(w^i_j)^2}$, falls beneath some user defined threshold---typically $5-10\%$ of $N$ \cite{SMCM_Doucet}. Note, unless $N_{eff}$ is large, this resulting discrete representation of the posterior is inherently coarse.

\subsubsection*{Perturbed-obs EnKF}: The ensemble Kalman filter with perturbed observations (summarized here 
following the work by Evensen ~\cite{ocdyn:evensen}) is a sequential data 
assimilation technique that evolves an ensemble of model states through time and 
performs Kalman filter style updates as new observations are incorporated.

Given an ensemble of $N_e$ model states at time $t_{j-1}$, each ensemble 
member is evolved according to \cref{model}. This forecast 
ensemble is used to generate a Gaussian estimate of the prior distribution at 
time $t_j$. We denote the forecast ensemble as $ \{\bx^i_{j,f} | 
i=1,...,N_e\}$. The forecast ensemble sample mean \(\mathbf{\bar{x}}_{j,f}\) 
and sample covariance \(\bP_{j,f}\) can be estimated as follows:
\begin{align}
\label{ENKF_forecast_ensemble_mean}
\bar{\bx}_{j,f} & = \frac{1}{N_e}\sum_{i=1}^{N_e} \bx^i_{j,f} \\
\label{ENKF_forecast_ensemble_covariance}
\bP_{j,f} & = \frac{1}{N_e - 1}\sum_{i=1}^{N_e} (\bx^i_{j,f} - 
\bar{\bx}_{j,f})(\bx^i_{j,f} - \bar{\bx}_{j,f})^T.
\end{align}

Observations are assumed to have the form $\bY_j = \bH\bx_j+\eta_j$, where $\bH$ is an observation matrix (typically a linearized  observation  operator) and observation errors $\eta_j$ are taken to be iid Gaussian random variables with mean 0 and known covariance $\bR$, 
i.e.~\(\eta_j \sim \mathcal{N}(\mathbf{0},\bR)\). We create an ensemble of $N_e$ 
perturbed observations with mean equal to \(\bY_j\) and covariance \(\bR\) 
according to $\bY^i_j = \bY_j + \epsilon^i_j$ 
where $\epsilon^i_j\sim \mathcal{N}(0,\bR)$. The covariance of 
the ensemble of perturbed observations is given by
\begin{align}
\label{ENKF_perturbed_observation_covariance}
\bR^e_j = \frac{1}{N_e-1} \sum_{i=1}^{N_e} \epsilon^i_j {\epsilon^{i}_j}^T .
\end{align}
The ensemble members are then updated according to
\begin{align}
\label{ENKF_update_eq}
\bx^i_{j,a} = \bx^i_{j,f} + \bP_{j,f} \bH^T ( \bH \bP_{j,f} \bH^T + \bR^e_j)^{-1} (\bY^i_j - 
\bH \bx^i_{j,f}) 
\end{align}
and the sample analysis mean and analysis covariance can be calculated as above yielding
\begin{align}
\bar{\bx}_{j,a} & = \frac{1}{N_e}\sum_{i=1}^{N_e} \bx^i_{j,a}, \quad \text{and} 
\quad \bP_{j,a} = \frac{1}{N_e - 1}\sum_{i=1}^{N_e} (\bx^i_{j,a} - 
\bar{\bx}_{j,a})(\bx^i_{j,a} - \bar{\bx}_{j,a})^T .
\end{align}
The analysis ensemble is used to generate a Gaussian approximation of the posterior distribution at time $t_j$. The analysis ensemble 
$\{\bx^i_{j,a}\}$ is then evolved to the next observation time by \cref{model} and used as the 
forecast ensemble for the next assimilation step. 

\subsection{Gaussian process emulators} \label{sec:gp}
The key approach in this work is, at time $t_j$, to learn about the mapping from an ``input space" (parameter space and/or state space at time $t_{j-1}$) to an ``output space" (state space) through the limited number of particle/ensemble samples. (Note, the weights of the particles do not inform us about this mapping directly.) To this end, we will employ a weakly stationary Gaussian process (GP) to model such an unknown relationship, $\bx_j \approx \hat f_j(\bx_{j-1}, \theta)$ or $\bx_j \approx \hat f_j(\theta)$. In the statistical computer models community, such modeling is typically referred to as statistical surrogates or GP emulators -- effectively statistical models of physical models. \cite{Sant:Will:Notz:2018,Rasm:Will:2006,Welc:Buck:etal:1992} provide excellent and broad overviews of this approach, but for the unfamiliar reader, we summarize the salient points here.

Consider $n_D$ {\em training} or {\em design} input values, $\bq^D \! =\{\bq_1,\ldots,\bq_{n_D}\}$, with each $\bq_k \in \mathbb{R}^r$, and a scalar output $y_k$ (e.g. ``output" may be one of the state variable values at time $t_j$, $y_k=x^k_j$) at each of these $n_D$ 
inputs, $\by^D=(y_1,\ldots,y_{n_D})^T$. We can model $\hat y \sim \MVN\big{(}m(\bq^D), \Sigma\big{)}$, a multivariate normal with $m(\cdot)$ a known mean trend and $\Sigma=\sigma^2 \gbR$, with variance $\sigma^2$. Here the correlation matrix $\gbR$ is computed by evaluating a chosen correlation function $c(\cdot,\cdot)$, e.g.~each element is given by $(\gbR)_{ij}=c(\bq_i,\bq_j)$. A Gaussian process emulator provides a prediction $\hat y (\bq^*)$ at an untried value of the input space $\bq^*$ as
\begin{align}
\nonumber \hat y (\bq^*)~=&~ h(\bq^*) \beta+\br^T(\bq^*) \gbR^{-1} (\by^D -h(\bq^D) \beta)+\delta \\
=&~ \hat f(\bq^*) \;. \label{gasp}
\end{align}
where $\br(\bq^*)=\big{(}c(\bq^*,\bq_1), \ldots, c(\bq^*,\bq_{n_D})\big{)}^T$. In other words, this gives the mean value of a Gaussian process at input $\bq^*$, where the process is conditioned to take on values of $\by^D$ at inputs $\bq^D$ if the uncorrelated noise term, $\delta$, is zero. Here $h$ is a set of basis functions (typically taken to be constant or linear), so $m(\bq)=h(\bq) \beta$ gives the overall trend based on the data, and the coefficient(s) are given by
\[
\beta = (h^T (\bq^D) \gbR^{-1} h(\bq^D) )^{-1} {h^T (\bq^D)} {\gbR^{-1}}{\by^D}.
\]
In these formulae, $\gbR$ is the $n_D \times n_D$ correlation matrix of the input design; often a power exponential or Mat{\'e}rn correlation kernel is assumed and ``fitting" an emulator amounts to finding the trend coefficients and correlation length scales that best represent the design pairs $\{\bq^D, \by^D\}$. We can also gain a sense of uncertainty induced by using the GP instead of the computer model simulation directly at $\bq^*$ by considering the standard prediction error
\[
s^2(\bq^*)=\sigma^2\Big{(} 1 - \br^T\gbR^{-1}\br + \frac{(1-\mathbf{1}^T\gbR^{-1}\br)^2}{\mathbf{1}^T\gbR^{-1}\mathbf{1}}\Big{)},
\]
where $\mathbf{1}$ is an $n$-vector of ones and $\sigma^2$ is the variance scaling of the process and found during the ``fitting" of the GP. Implementations of GP emulators are available: in \textsc{Matlab} one can use the function~\texttt{fitrgp()}\footnote{Introduced in version R2015b, see \texttt{https://mathworks.com/help/stats/fitrgp.html}}, or in~\textsc{R} the \texttt{Robustgasp()} package \cite{Gu19}.

The Parallel Partial Emulator (PPE) generalizes the standard emulator construction presented in \cref{gasp} for scalar outputs, to an emulator for vector-valued outputs \cite{gu:berg:2016}. Consider then a set of $N$ model design inputs and $s$-dimensional responses $\{\bq^D, \bm{Y}^D\}$. 
$\bm{Y}^D$ is now an $N\times s$ matrix. PPE allows each output component to have a unique mean $m_j(\bq)=h(\bq) \psi_j$ and variance $\sigma^2_j $ ($j=1,\ldots,s$), but assumes a shared correlation structure and correlation parameters among all locations. Equations for predictive mean and standard error are nearly identical to \cref{gasp}, but are $s-$dimensional. 
We mention that the means and variances of the individual Gaussian processes inherit some measure of (spatial) correlation that is present in the physical system, even though no explicit assumption is made about spatial relationships. 

\section{Methodology}
\label{sec:method}
This section constructs approximations of the Particle Filter that employ only a relatively small number of model runs. The model runs are used as design-response pairs in a Gaussian Process emulator; a large number of samples from the GP emulator are then treated as particles in a PF. Several algorithms are presented here, as the practical options for emulator design and response variables depend on the parameter and state dimension. 

The following \cref{sec:emupf} introduces a naive but straightforward blending of the GP emulator and PF, which is then employed as a springboard to introduce multiple refinements.

We employ subscripts for time indices and superscripts for particle indices, and additionally use superscripts in parentheses to denote components of a vector.

\subsection{The Emulator Particle Filter: Emu-PF}
\label{sec:emupf}
We construct an emulator for the map \cref{model} from time $t_{j-1}\rightarrow t_j$. Then, we use the emulator output in a PF as if it were samples from the prior distribution in \cref{bayes}.

At time $t_{j-1}$ suppose we have evenly weighted parameter estimates and state estimates $\theta^i_{j-1}\in\R^{p}$ and $\bx^i_{j-1}\in\R^{n}$, $i$ from $1$ to $n_D$. Suppose also we have a large ensemble of parameter estimates $\Theta^i_{j-1}$ and corresponding state estimates $\bX^i_{j-1}$ with weights $w^i_{j-1}$, $i$ from $1$ to $N_F$. Then follow the following sequence:
\begin{enumerate}
\item[\bf Forecast:] Employ the numerical model \cref{model},
    \[ \bx_j^i = \varphi(\bx_{j-1}^i,\theta_{j-1}^i) \;, \]
    Label the $l$th component of $\bx^i_j$ by $\bx^{i,l}_j.$ Set $\theta^i_{j} = \theta^i_{j-1}$ and $\Theta^i_{j} = \Theta^i_{j-1}$ (or employ a parameter model).
\item[\bf Emulate:] For each $l$ from $1$ to the state dimension $n$,
    \begin{enumerate}[label=E\arabic*.,ref=E\arabic*]
    \item\label{emuSetI} Set the emulator design variables (inputs) to be the state and parameter estimates at which we employed the model,
    \[
    \bq^D = \{x^i_{j-1},\,\theta^i_{j-1}\}_{i=1}^{n_D} \;.
    \]
    \item\label{emuSetO} Set the response variables (outputs) to be the $l$th component of the state variable from the model output,
    \[
    \by^D = \{\bx^{i,l}_j\}_{i=1}^{n_D} \;,
    \]
    and fit the emulator with the design-response pairs.
    \item\label{emuEval} Evaluate the emulator at each of the fine state and parameter values: set
    \[
    \bq^* = \{\bX^i_{j-1},\,\Theta^i_{j-1}\}_{i=1}^{n_D} \;, 
    \]
    obtain $\hat y(\bq^*) = \hat f(\bq^*)$ from \cref{gasp}, and save each scalar emulator output as the $l$th component of $\bX^i_j$, $\bX^{i,l}_j = \{\hat y(\bq^*)\}_i$, for each $i$ from $1$ to $N_F$.
    \end{enumerate}
\item[\bf Assimilate:] Treat $\{\Theta^i_j, \bX^i_j \}$ as samples from the prior and perform a Data Assimilation scheme in the high-dimensional $N_F$ space. For example, in a Particle Filter, employ \cref{w} to judge the emulator outputs,
    \[ w_j^i = \frac{
    \exp(-\frac{1}{2}\left(\by_j-h(\bX_j^i)\right)^T \bR^{-1} \left(\by_j-h(\bX_j^i)\right)
    }{
    \sum_{k=1}^N\exp(-\frac{1}{2}\left(\by_j-h(\bX_j^k)\right)^T \bR^{-1} \left(\by_j-h(\bX_j^k)\right)
    } \;,\]
    assuming the observation errors are Gaussian with covariance $\bR$, and the observation operator is $h()$. Calculate the effective sample size $N_{eff}$, defined below \cref{w}, and resample $\{\Theta^i_j, \bX^i_j \}$ if needed.
    \item[\bf Subsample:] Use a resampling algorithm to sample $n_D$ times from the weighted pairs $\left( \{\Theta^i_j, \bX^i_j \},\, w_j^i\right)$.
\end{enumerate}
In the above we use $\{.\}_{i=1}^{N}$ to indicate when an operation can be vectorized by concatenating together ensemble members as columns in a matrix. 

\begin{figure}
    \centering
\definecolor{Mcol}{RGB}{217,  95,   2}
\definecolor{Mtxt}{RGB}{165,  71,   2} 
\definecolor{mcol}{RGB}{0,  123,   123} 
\definecolor{mtxt}{RGB}{0,  100,   100} 
\definecolor{excol}{RGB}{76,  0,   153}  
\definecolor{extxt}{RGB}{80,  0,   163}  

\def\fwid{\linewidth*.4} 
\def\nDist{3em} 
\def\bHgt{1em} 

\tikzstyle{decision} = [rectangle, draw, fill=Mtxt!20, 
     text badly centered, inner sep=4pt]
\tikzstyle{cloud} = [draw, ellipse,fill=Mtxt!20, node distance=4em,
    minimum height=\bHgt]
\tikzstyle{block} = [rectangle, draw, fill=mcol!20, 
     text centered, rounded corners, minimum height=\bHgt]
\tikzstyle{line} = [draw, thick, -latex']
\small{
\begin{tikzpicture}[node distance = \nDist, auto]
    \node [decision] (dim) {Ensemble size:};
    \node [below of=dim] (mPos) {};
    \node [cloud, right=\fwid/4 of dim] (nd) {$n_D$};
    \node [cloud, right=0.9*\fwid of nd] (nf) {$N_F$};
    
    \node [block, below of=nd](ncrse) {Ensemble at time $t_{k-1}$};
    \node [block, below of=nf](nfne) {Ensemble at time $t_{k-1}$};
    
    \node [block, below of=ncrse] (crsemod) {Model forecast to time $t_{k}$};
    \node [block, below of=crsemod,fill=extxt!20] (emuGen) {Construct $t_{k-1}\rightarrow t_{k}$ emulator};
    
    \node [block, below=(2*\nDist-1.6*\bHgt) of nfne,fill=extxt!20] (emuUse) {Evaluate emulator for large ensemble};
    
    \node [block, below of=emuUse] (pf) {Run Data Assimilation scheme};
    
    \node [block, below of=pf] (fneOut) {Ensemble at time $t_k$};
    \node [block, below=(2*\nDist-1.5*\bHgt) of emuGen] (crseOut) {Subsample ensemble at time $t_{k}$};
    

    \path [line] (ncrse) ->  (crsemod);
    \path [line] (crsemod) ->  (emuGen);
    \path [line] (emuGen) ->  (emuUse);
    \path [line] (nfne) ->  (emuUse);
    \path [line] (emuUse) ->  (pf);
    \path [line] (pf) ->  (fneOut);
    \path [line] (fneOut) ->  (crseOut);

\end{tikzpicture}
}
    \caption{Overview of the novel synthesis of Gaussian process emulators with Data Assimilation methods.}
    \label{fig:flow}
\end{figure}
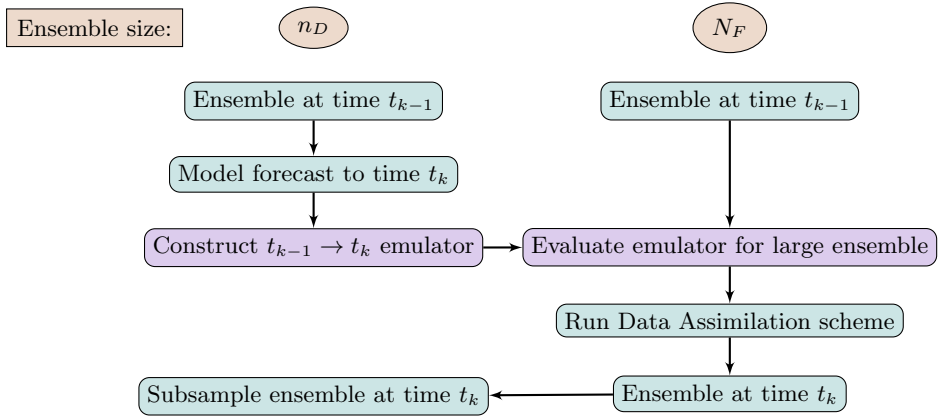

\begin{figure}
    \centering
       \begin{subfigure}[b]{\textwidth}
    \includegraphics{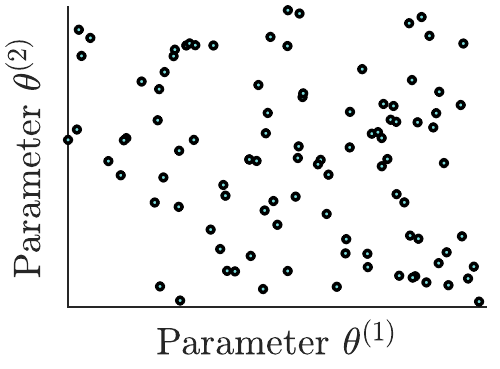} \includegraphics{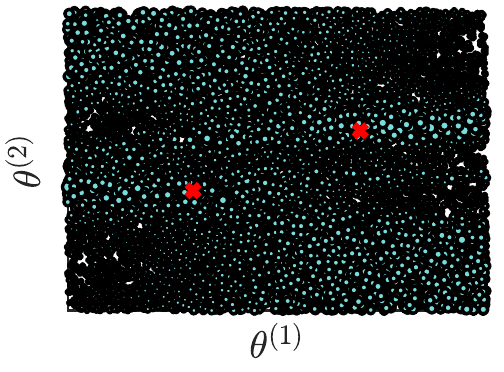}
    \caption{Parameter design variables (left) and fine weighted ensemble (right) before assimilation. Red crosses show modes of the true parameter distribution.}
    \label{fig:subABt1}
    \end{subfigure}
           \begin{subfigure}[b]{\textwidth}
    \includegraphics{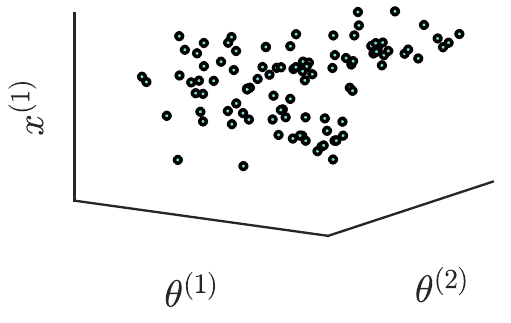}  \includegraphics{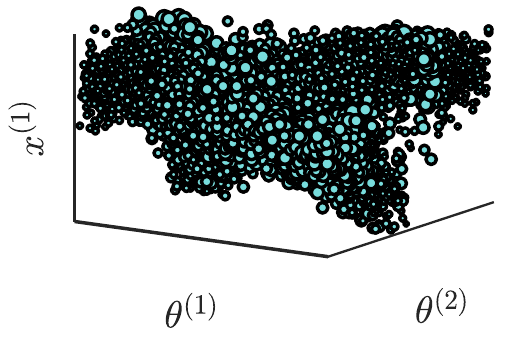}
    \caption{Left: state variable $x^{(1)}$ after model integration, viewed as a function of two parameters. This is used as input to train the emulator. Right: emulator output at $10,000$ samples. This output has been weighted by a measure of distance from the observations (the highest-weight particles, with larger dot sizes, are concealed within the cloud of samples).  }
    \label{fig:subX}
    \end{subfigure}
               \begin{subfigure}[b]{\textwidth}
    \includegraphics{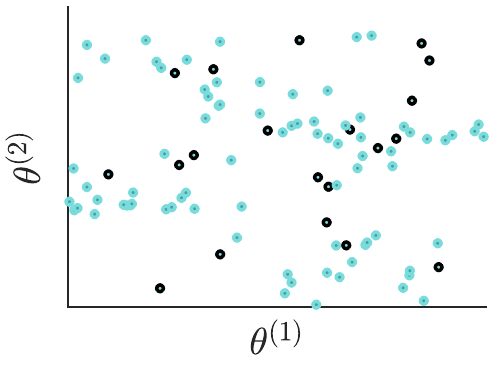}  \includegraphics{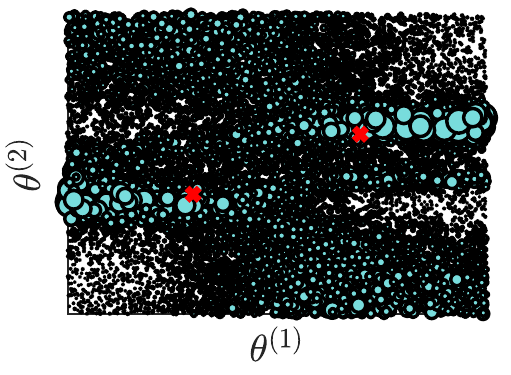}
    \caption{Parameter estimates after assimilation. Right: weights in the fine ensemble, shown by dot size, reveal a highly nonGaussian posterior that is shifting towards the true parameter values (red crosses). Left: the brighter dots have been sampled from the fine ensemble; the darker dots are unchanged from Figure~\textsc{(a)}. \Cref{sec:algorithms} explains the benefits of keeping some of the design variables fixed. }
    \label{fig:subABt2}
    \end{subfigure}
    \caption{Visualisation of the internal Emu-PF mechanisms over one assimilation step. Left column shows components of dimension $n_D=100.$ Right column shows components of dimension $N_F=10,000$. \textsc{(a)}: parameter ensembles at time $t_j$. \textsc{(b)}: distribution of one state variable as a function of parameters. \textsc{(c)}: parameter ensembles at time $t_{j+1}$.  Full details for this $8$ state, $2$ parameter experiment are given in \cref{sec:numerics}.} 
    \label{fig:emuInternal}
\end{figure}

\begin{figure}
    \centering
    \includegraphics{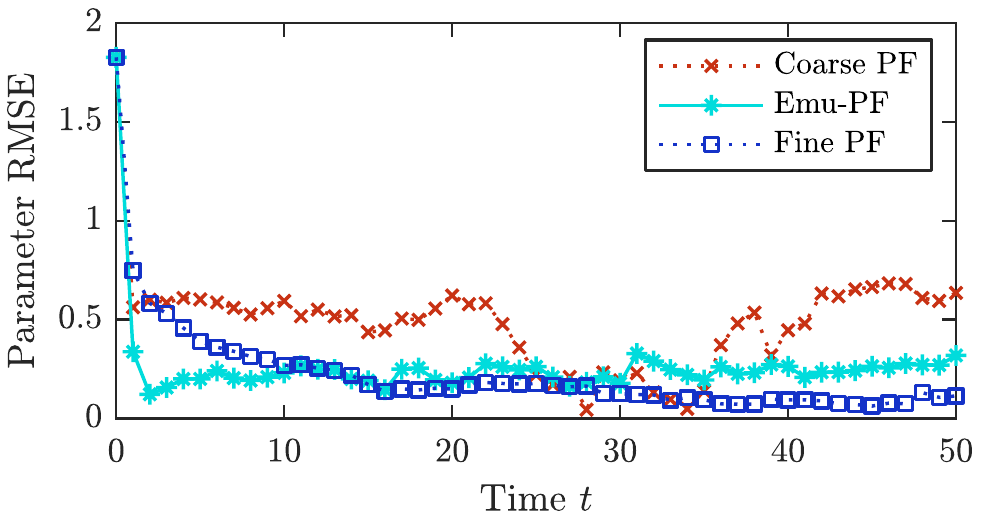}\\
    \includegraphics[]{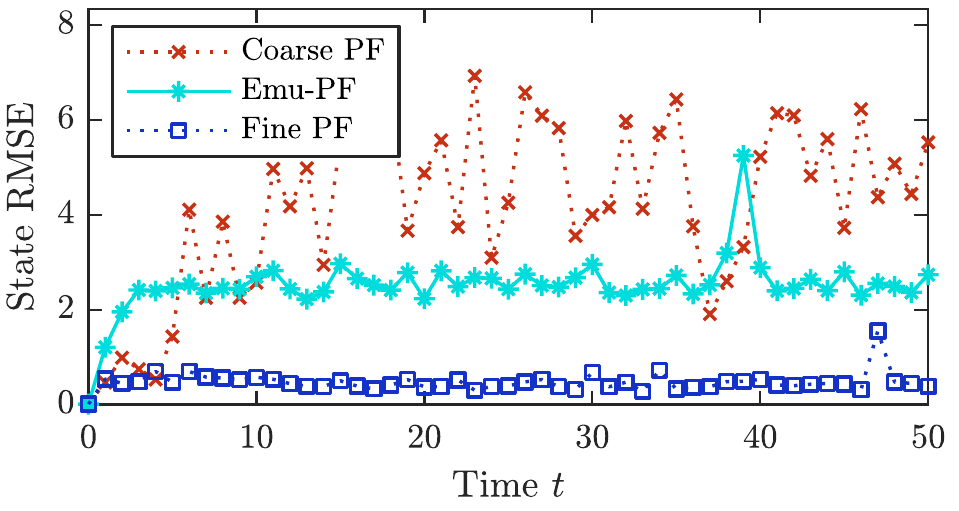}
    \caption{Long term error statistics for the implementation of Emu-PF from \cref{fig:emuInternal}, compared to: a ``coarse'' PF that employs $n_D=100$ model runs as done in the Emu-PF, and a ``fine'' PF that employs $N_F=10,000$ model runs, equal to the number of samples in the Emu-PF emulator. For this implementation, performance of Emu-PF is markedly better than the coarse PF but not to the standard of the fine PF. Modifications to Emu-PF that empower it to compete with the fine PF, despite employing only $1\%$ as many model runs, are discussed in \cref{sec:Epf,sec:ppepf,sec:pcapf} and displayed in \cref{sec:numerics}.}
    \label{fig:emuSampleRMS}
\end{figure}

\Cref{fig:flow} shows a schematic for this class of surrogate DA methods. The key steps in this schematic are displayed for an Emu-PF (with implementation details delayed until \cref{sec:numerics}) in \cref{fig:emuInternal}. The long-time error statistics for this Emu-PF are compared to Particle Filters employing $n_D$ particles and $N_F$ particles in \cref{fig:emuSampleRMS}. 

Two shortcomings of GP emulators motivate improvements in the above algorithm. First, it is notoriously challenging to fit emulators with a high-dimensional input space. Yet the surrogate DA method employs high-dimensional inputs to the emulator, as the parameter and state vectors are combined and used as design variables. Some recent works \cite{betancourt2020gaussian,liu2017dimension} offer approaches for dimension reduction for statistical emulators that require either significant prior knowledge of the variability of the input space or a significant amount of data to characterize that variability well. For sequential DA, as a matter of course we have this prior knowledge available, but the flavor of appropriate dimension reduction will be problem specific. In particular, DA schemes that employ localization may favor a dimension reduction approach that is local as opposed to a global dimension reduction. We will explore variations of each.
Secondly high-dimensional response variables are avoided by looping through the entire state vector, one dimension at a time; but this is a potentially slow and expensive procedure. There are multiple recent approaches to emulating high-dimensional output and we will explore a variation of our algorithm that utilizes one of those approaches, namely partial parallel emulation \cite{gu:berg:2016}.

We now introduce various variations on the Emu-PF. Each either reduces the dimension of the design variables or improves the efficiency of sampling from the emulator.

\subsection{Variant: include only some state values in the emulator input} 
\label{sec:Epf}
This modification implements a straightforward localization for the emulator inputs. Modify the emulation step to include only state values near the response variable. For each $l$ from $1$ to the state dimension $n$,
\begin{enumerate}[label=Es\arabic*.,ref=Es\arabic*]
\item\label{E.emueval} Choose some integer $\Gamma$. Set the emulator design inputs to be the parameter estimates at which we employed the model, and a slice of the state inputs,
    \[
    \bq^D = \left\{\left(\theta^i_{j-1},\,\bx^{i,\,(l-\Gamma:l+\Gamma)}_{j-1}\right)\right\}_{i=1}^{n_D} \;.
    \]
\item As \cref{emuSetO}.
\item Evaluate the emulator at each of the fine parameter values and corresponding state estimates: set
    \[
    \bq^* = \left\{\left(\Theta^i_{j-1},\,\bX^{i,\,(l-\Gamma:l+\Gamma)}\right)\right\}_{i=1}^{n_D} \;, 
    \]
    otherwise as \cref{emuEval}.
\end{enumerate}
In 2-d or 3-d space, instead choose a localization distance parameter $\Gamma\ge0$ and include every grid point within radius $\Gamma$ of $\bx^{.,(l)}_j$ in the design input step \cref{E.emueval}.

One extremely simple implementation of this variation on the Emu-PF is to set $\Gamma=-1$; that is, to include no state variables at all in the emulation. This implementation is justified if the distribution $x|\theta\approx g(\theta)+$noise, for a smooth function $g$. Equivalently, the distribution $x|\theta$ should be roughly unimodal. This condition is frequently satisfied in practice \cite{Morzfeld19}, and the resulting algorithm is fast but still readily capable of filtering nonGaussian marginal distributions for $\theta$. 

If this variant of the Emu-PF is employed, we refer to it by the value of $\Gamma$ chosen; all DA methods are benchmarked against the Emu-PF with $\Gamma=-1$ in \cref{fig:m8,fig:m2,fig:os5m8,fig:os5m4}, and \cref{fig:os5m8,fig:os5m4} test the Emu-PF with $\Gamma=1$.

\subsection{Variant: compute emulator outputs in parallel with \texttt{ppgasp}}
\label{sec:ppepf}
Use Partial Parallel Estimation (described in \cref{sec:gp}) to compute all states at once.
\begin{enumerate}[label=Ep\arabic*.,ref=Ep\arabic*]
\item \label{emuPi} Set the emulator design inputs to be the parameter estimates at which we employed the model,
    \[
    \bq^D = \{\theta^i_{j-1}\}_{i=1}^{n_D} \;.
    \]
\item \label{emuPo} Set the response variables to be the model output,
    \[
    \by^D = \{\bx^i_j\}_{i=1}^{n_D} \;.
    \]
\item \label{emuPe} Evaluate the emulator at each of the fine parameter values: set
    \[
    \bq^* = \{\Theta^i_{j-1}\}_{i=1}^{n_D} \;, 
    \]
    obtain $\hat y(\bq^*) = \hat f(\bq^*)$ from \cref{gasp}, and save each column of emulator output as $\bX^i_j = \{\hat y(\bq^*)\}_i$ for each $i$ from $1$ to $N_F$.
\end{enumerate}
The above implementation avoids the for-loop present in \cref{sec:emupf,sec:Epf}, but as written it is incompatible with the $\Gamma$-localization discussed in \cref{sec:Epf}, as the PPE emulates all state variables simultaneously. We discuss simultaneous parallelization and localization in \cref{sec:discussion}. A more radically localized, efficient Emu-PF for state estimation is discussed in \cref{sec:Xpf}.

\subsection{Variant: perform a global dimension reduction before using emulator inputs}
\label{sec:pcapf}
Employ a data-based dimension reduction algorithm (e.g. PCA, DMD, diffusion maps, UMAP, \dots) on the state variables going into the emulation mapping of $\bx_{j+1}=\varphi(\bx_{j};\theta_j)$. This approach is not generically used to emulate high dimensional parameter inputs because it's often unclear how to represent the variability of parameters, but in the sequential DA case there is an obvious candidate---the vector of state variables.

As a clear example, in the remainder of the section and in numerical examples we employ PCA. That is, we have an approximation from the fine sampled posterior at the $j$th time step, $\left( \{\Theta^i_{j},\bX^i_{j}\},\, w_j^i\right)_{i=1}^{N_F}$. Let $X=X_{data}-\bar X_{data}\mathbf{1}_{N_F}$ be the $n\times N_F$ matrix where the $i$th column of $X_{data}$ is $\bX^i_j$, $\bar X_{data} =\frac{1}{N_F} \sum_{i=1}^{N_F} \bX^i_{j}$, and $\mathbf{1}_{N_F}$ is a row vector consisting of $N_F$ ones. Then $A=XX^T$ is a covariance matrix representative of the variance in $X$. A singular value decomposition of $A$ produces $A=V \Lambda V^T$ where $\Lambda$ is a unitary matrix containing ordered singular values, the columns of $V$ contain the corresponding singular vectors, and $V^T=V^{-1}$ as $A$ is symmetric. Truncate $\Lambda$ and $V$ to keep only the largest $r<n$ singular values; label the truncated matrices $\tilde \Lambda$, now $r\times r$, and $\tilde V$, now $J\times r$. Note $A \approx \tilde V \tilde \Lambda \tilde V^T$. 
Now let $Y= \tilde V^T X$. 
Effectively $Y$ is a matrix of weights to multiply the principal components vectors (columns of $\tilde V$) to recover the original data $X$. 

In Emu-PF schemes employing PCA, we use the weights $Y$ as input variables for emulation in \cref{emuSetI}. The response variables are unchanged in \cref{emuSetO}, but when evaluating the emulator at fine samples in \cref{emuEval} we replace $\bX_{j-1}^i$ with $\tilde V^T \bX_{j-1}^i$. 

We discover a fast, flexible and powerful Emu-PF algorithm by combining global dimension reduction of inputs (by PCA in our experiments) and fast emulator outputs (with PPE, described in \cref{sec:ppepf}); this algorithm is employed in Experiments Two and Four of \cref{sec:numerics}. 

\subsection{Variation: localize the emulator by ``slicing and stacking'' the emulator inputs}
\label{sec:Xpf}
This variation on the Emu-PF involves a radical rethinking of the emulator state inputs; for that reason we suppress parameter dependence and consider state estimation only. Assume that the physical law governing state evolution is the same for each component of the state vector; then a single model run, initialized at $x_{j-1}^i$ and producing $x_j^i\in\R^n$, provides $n$ samples of that physical law. The following algorithm exploits this rich data by configuring the emulator design inputs as $n\times n_D$ samples, rather than $n_D$ samples. 

We suppose that some localized slice of state variables at time $j-1$, within distance $\Gamma$ of state variable $l$, is sufficient to predict the $l$th state variable at time $j$. The following procedure learns a $\R^{2\Gamma+1}\rightarrow\R$ map for the state update.
\begin{enumerate}[label=Er\arabic*.,ref=Er\arabic*]
\item Choose some integer $\Gamma\ge0$. The design inputs $\bq^D$ are to be a {$(2\Gamma+1)\times (n\times n_D)$}~array, with the $q$-th row of that array given by
    \[
    \left\{\bx^{i,\,l-\Gamma:l+\Gamma}_{j-1}\right\}_{i=1}^{n_D} \;,
    \]
    where $i=\ceil(q/n)$ and $l\;=\mod(q,n).$ 
    \item Set the response variables $\by^D$ to be the corresponding $n\times n_D$-vector of state variables, with the $q$th entry
    \[
    \{\bx^{i,l}_j\}_{i=1}^{n_D} \;,
    \]
    \item Evaluate the emulator at each of the state estimates: set
    \[
    \bq^* = \left\{\bX^{i,\,l-\Gamma:l+\Gamma}\right\}_{i=1}^{n_D} \;, 
    \]
    otherwise as \cref{emuEval}.
\end{enumerate}
This approach entails a radical reduction in the dimension of emulator inputs \emph{and} outputs. Due to the unusual ``slicing'' of the emulator input to obtain rich training data, we refer to it as the ``sliced Emu-PF.'' We test it on a state estimation problem in \cref{fig:m4s}.

 \section{Numerical experiments and results} 
 \label{sec:numerics}
 
 We consider a joint state-parameter estimation problem from \cite{SantitissadeekornJones15}. The state $\bx_j$ is generated by integrating from time $ t_{j-1}$ to $ t_j$ the system of ordinary differential equations introduced in \cite{Lo96},
\begin{align} \label{l96}
    \dot x^{(l)} =& \left(x^{(l+1)}-x^{(l-2)}\right)x^{(l-1)} -x^{(l)} +F^{(l)}\;,
\end{align}
commonly called the Lorenz-96 system. Superscripts in parentheses denote components of a vector, $l$ ranges from $1$ to $n$, and the forcing depends on two parameters
\begin{align}
    \label{forcing}
 F^{(l)} = 8 + \theta^{(1)}\sin\left(\frac{2\pi l}{n\theta^{(2)}}\right) \,. 
 \end{align}
 We will compare the surrogate DA algorithms to Particle and Ensemble Kalman Filters. Our goal is to obtain performance similar to that of a Particle Filter that employs a large number of particles, $N_F$, but only allowing $n_D \ll N_F$ model runs in our scheme. In order to quantify the benefits, and drawbacks, of our approach, we will include the following algorithms for comparison:
 \begin{itemize}
     \item A ``fine PF'' employing $N_F$ particles,
     \item A ``coarse PF'' employing $n_D$ particles,
     \item An EnKF employing $n_D$ particles.
 \end{itemize}
While several of our results feature implementations of the Emu-PF that compete with, or exceed the performance of, the fine PF, it is important to remember that our original goal was to attain performance somewhere between the coarse and fine PF. Exact implementation details for all DA methods are given in \cref{sec:algorithms}. We will also briefly discuss better implementations of the Particle Filter. For clarity in the results we do not employ any of these advanced filtering methods, but emphasize that our approach does not conflict with the usage of them. 
 
 It is standard in the atmospheric forecasting community to employ \cref{l96} with dimension $n=40$, and to compute and subsequently discard a ``burn-in'' period of at least a thousand assimilation steps. Our benchmark fine PF is incapable of resolving the $n=40$ case without extensive modifications that, if also implemented in an Emu-PF, can make it difficult to be sure what the contributions of the emulator are. Additionally, good filter performance during the first twenty assimilation steps are crucial for parameter estimation (assuming an initially uninformative prior on the parameters). For these reasons we choose model dimension $n=8$, analysed in \cite{Orrell03}, and include the filter performance over the initial assimilation steps.
 
 Over all experiments, a vast quantity of information is computed. We will summarize this information with the Root Mean-Squared-Error (RMSE) and the sample variance. For parameter estimates the posterior distribution is multimodal (see \cref{sec:exp}); when calculating RMSE or variances of parameter estimates, we first apply absolute values to reduce the number of modes.
 
 \subsection{Implementation Details}
 \label{sec:algorithms}
 \emph{Particle Filters} all employ the merging particle filter of \cite{Nakano07}, with the recommended values $a_1 = 3/4$; $a_2=(\sqrt{13}+1)/8$; $a_3 = -(\sqrt{13}-1)/8$. Additionally PFs employ the parameter model from \cite{Liu01} which, at each observation time, jitters all particles randomly by adding noise generated with standard deviation $\alpha=0.01$, then draws all particles slightly towards the particle mean, preserving $\beta=0.99$ of the variation among particles. 
 
 \emph{The EnKF} employs multiplicative covariance inflation of $1.02$. That is, when the sample forecast covariance is calculated in \cref{ENKF_forecast_ensemble_covariance}, it is multiplied by $1.02$ before it is used in \cref{ENKF_update_eq}. Covariance inflation is a common remedy to the problem of a slightly under-dispersive ensemble in the EnKF.
 
 \emph{The Emu-PF} algorithms divide the $n_D$ particles that are used in model runs into two groups. The first group {is} sampled from the fine posterior after every assimilation step, as described earlier in the paper (\cref{fig:flow}, for example). The second group is not sampled from the posterior, and remains fixed over the assimilation steps. We fix this second group, comprising $20$ of the $100$ design variables, so that the emulator can evaluate inputs at a wide range of $\theta$ even if the subsampled group has narrow support. \Cref{fig:subABt2} (left) shows the first group ($80$ bright dots) and second group ($20$ dark dots). 
 
 A modern implementation of the particle filter to a high-dimensional filtering problem should involve intensive modifications to mitigate the curse of dimensionality. Successful innovations include proposal densities \cite{van19}, mixtures \cite{Crisan15}, and dimension reduction strategies including the classic Rao-Blackwellized PF or recent localized PFs \cite{Poterjoy16,Potthast19}. For clarity of exposition we employ none of the above modifications to the PF, but note that the Emu-PF is compatible with, and different in approach to, all of them (several examples are discussed in \cref{sec:discussion}).

\subsection{Experiment Details}
\label{sec:exp}
For all experiments we assimilate data at $1000$ observation times with time step of $0.05$ between them. Model and truth are integrated between these observation times with five steps of a fourth order Runge-Kutta scheme. At each of these integration steps the true value of $\theta^{(1)},\,\theta^{(2)}$ is drawn from a Gaussian with mean $(2,1)^T$ and variance $0.01\mathbf{I}_2$. All DA schemes use fixed parameter estimates between assimilation steps. The discrepancy between fixed parameters in model updates for DA schemes, and varying parameters by drawing them from a distribution for data introduces a simple form of model error. We initialize state ensembles at $t=0$ with a tight spread of variance $0.01\mathbf{I}_8$ around the true initial condition, which is generated randomly. By contrast the parameter ensembles are initially uninformative, being drawn from a uniform distribution on the square $(-5,5)\times(-5,5)$. The symmetry of \cref{forcing} ensures that the posterior distribution in the parameters is always at least bimodal, as the forcing $F^{(l)}$ is identical at $+\theta,\,-\theta$; but, also by symmetry, we can calculate reasonable RMSE and variance statistics for parameters by taking the absolute value of parameter estimates. All schemes employ $n_D=100$ and $N_F=10,000$. 

Interpret plots of the DA schemes with the following: if two different initial conditions for \cref{l96,forcing} are integrated for a long time, the mean distance between the two trajectories will be around $5$. Any DA method attaining a state RMSE value near $5$ is no different to employing no assimilation. However smaller RMSE is not necessarily optimal; each DA scheme is trying to estimate the posterior, which is unknown. Generally we will compare methods to results from the fine PF.

We also present tables with summary statistics for each experiment. These tables present the mean RMSE and the median sample variance over the final $50\%$ of assimilation steps, recorded separately for parameters and for states. We compute median variance as the mean variance is dominated by large variance terms in a few of the state variables. Generally the sample variances will appear to suggest methods are under-dispersive; but the EnKF performs better estimating the bimodal parameter distribution if it is under-dispersive than otherwise (explained further in the discussion of Experiment One).

We vary two quantities between experiments; the dimension $m$ of the observations, and the accuracy of the observations. We will consider $m=2,\,4,\,8$ evenly spaced observations. The observation accuracy is measured by the scalar $\sigma_o$, which controls the observation error covariance matrix $\bR$ from \cref{sec:da} according to $\bR = \sigma_o^2 \mathbf{I}_m$.
More difficult experiments are obtained by reducing $m$ and/or $\sigma_o$. Fewer observations at each observation time lead to a more uncertain posterior, which is difficult for the Emu-PF algorithms to represent with the low number of design variables $n_D$. Accurate observations, that is smaller values of $\sigma_o$, are difficult for Particle Filters in general.

\subsection*{Experiment One} 
\begin{figure}
     \centering
     \includegraphics{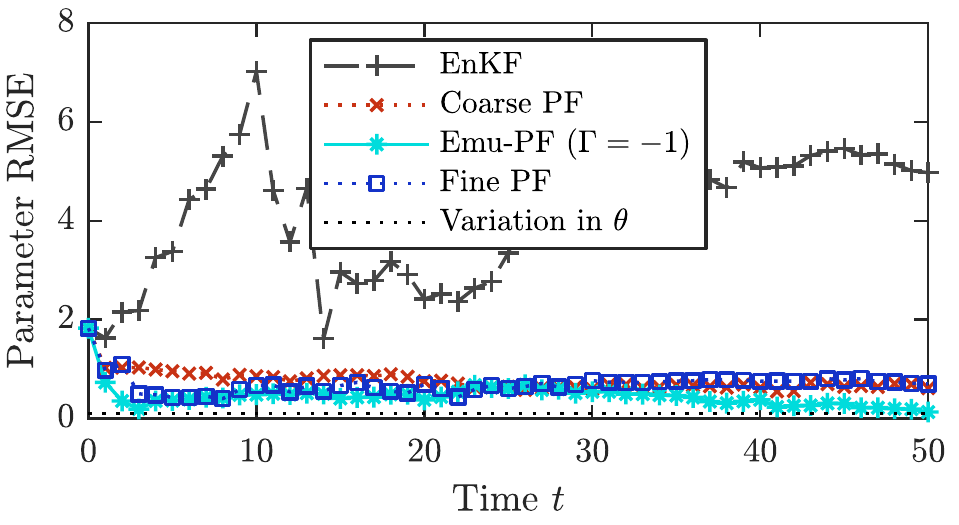}
     \includegraphics{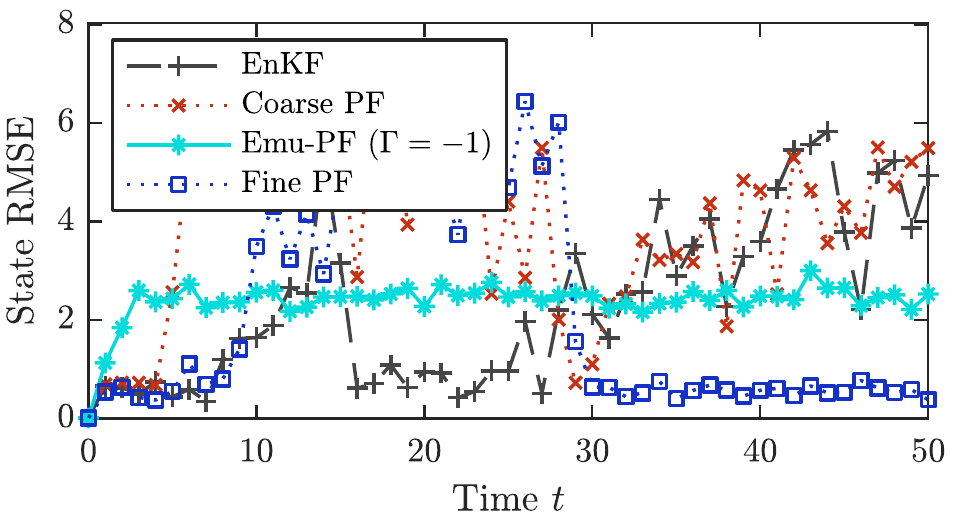}
     \caption{Error statistics for Experiment One, $m=8$ observations at each observation time, of accuracy $\sigma_0=1$. In this (and every) plot, only every $20$th data point is shown. For this mildly difficult filtering problem, we observe that the $\Gamma=-1$ implementation of \cref{sec:Epf}, that uses no state variables at all as emulator inputs, is stable and reasonably accurate.}
     \label{fig:m8}
 \end{figure}

We begin by presenting statistics for a fully observed ($m=8$) system with observation accuracy $\sigma_o=1.$ in \cref{fig:m8,table:m8}. The Emu-PF with $\Gamma=-1$ outperforms the, equivalent in number of model runs, coarse PF. The fine PF does not appear to estimate the state variables well in \cref{fig:m8}; \cref{fig:emuSampleRMS} is another run with the same setup in which the fine PF is clearly distinct from the coarse. \Cref{table:m8} shows that the Emu-PF with $\Gamma=-1$ has smaller variance in parameter estimates than the fine PF, but has larger variance, roughly equal to RMSE, in state variables. 

The EnKF performs poorly in \cref{fig:m8}, with large errors compared to other schemes in both the parameters and state variables; we now discuss why. Poor performance is expected, as the distribution of the parameters is bimodal and the EnKF relies on unimodal approximations. In practice we observed this poor behaviour from the EnKF only in about half of all experiments; in the remainder of experiments, the EnKF parameter ensemble tends to shrink over assimilation steps. Once the parameter ensemble has collapsed sufficiently, the peaks of the posterior---visible in \cref{fig:subABt2}, right---are no longer both contained in the span of the ensemble, and the ensemble can move close to one or another peak in subsequent assimilation steps. \Cref{table:m8} shows that the Emu-PF tends to overtrain on parameter estimates, but maintains a variance in state variables roughly equal to the corresponding RMSE.

\begin{table}                                                      
\centering                                                         
\begin{tabular}{ccccc}                                             
 & RMSE ($\theta$) & Var ($\theta$) & RMSE ($\bx$) & Var ($\bx$) \\
Fine PF & 0.74 & 0.0048 & 1.2 & 0.12 \\                            
Coarse PF & 0.65 & 0.0014 & 3.7 & 0.02 \\                          
EnKF & 4.8 & 0.024 & 3.4 & 0.039 \\                                
Emu-PF ($\Gamma=-1$) & 0.39 & 0.00028 & 2.4 & 2.9 \\               
\end{tabular}                                                      
\caption{Summary statistics for Experiment One.}                   
\label{table:m8}                                                   
\end{table}  

\subsection*{Experiment Two: sparse observations} 
\begin{figure}
     \centering
     \includegraphics{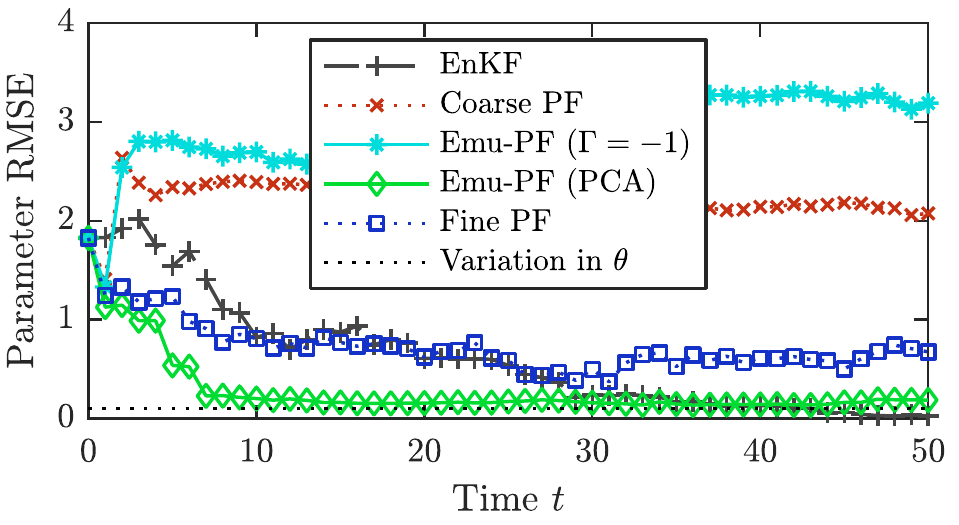}
     \includegraphics{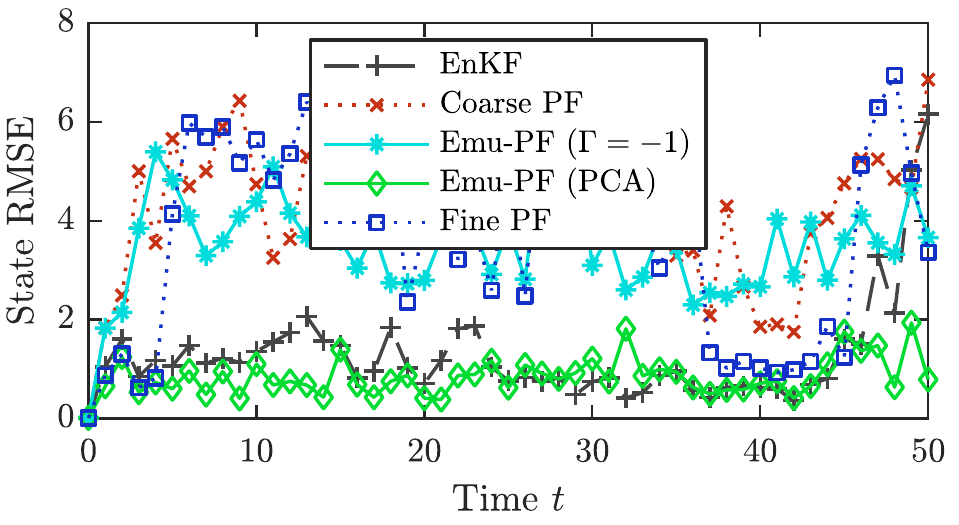}
     \caption{Error statistics for Experiment Two, $m=2$ observations at each observation time, of accuracy $\sigma_0=1$. The $\Gamma=-1$ Emu-PF and fine PF both under-perform compared to their mean behaviour; the Emu-PF employing PCA is stable and accurate.}
     \label{fig:m2}
 \end{figure}

Consider the more challenging setup of $m=2$ evenly spaced observations with the same accuracy $\sigma_o=1.$ from the previous experiment. An implementation of the Emu-PF employing both, the PCA dimension reduction from \cref{sec:pcapf} (to four variables), and PPE from \cref{sec:ppepf} to compute all emulator outputs simultaneously, is tested in \cref{fig:m2,table:m2}. The Emu-PF implementation with PCA not only stably estimates states and parameters under a difficult filtering problem, but out-competes both the EnKF and the fine PF (which employs $100\times$ as many model runs). This good performance from the Emu-PF with PCA far outstrips our original goal, that was just to \emph{replicate} the performance of the fine PF. 
 
 \begin{table}                                              
\centering                                                      
\begin{tabular}{ccccc}                                          
 & RMSE ($\theta$) & Var ($\theta$) & RMSE ($\bx$) & Var ($\bx$) \\
Fine PF & 0.57 & 0.0033 & 3.6 & 0.21 \\                         
Coarse PF & 2.2 & 0.0026 & 4.1 & 0.0046 \\                      
EnKF & 0.17 & 0.0018 & 1.2 & 0.25 \\                            
Emu-PF ($\Gamma=-1$) & 3.2 & 0.00068 & 3.4 & 3 \\               
Emu-PF (PCA) & 0.16 & 0.00099 & 1 & 0.25                     
\end{tabular}                                                   
\caption{Summary statistics for Experiment Two.}                
\label{table:m2}                                                
\end{table} 
 
 \subsection*{Experiment Three: accurate observations} 
\begin{figure}
     \centering
     \includegraphics{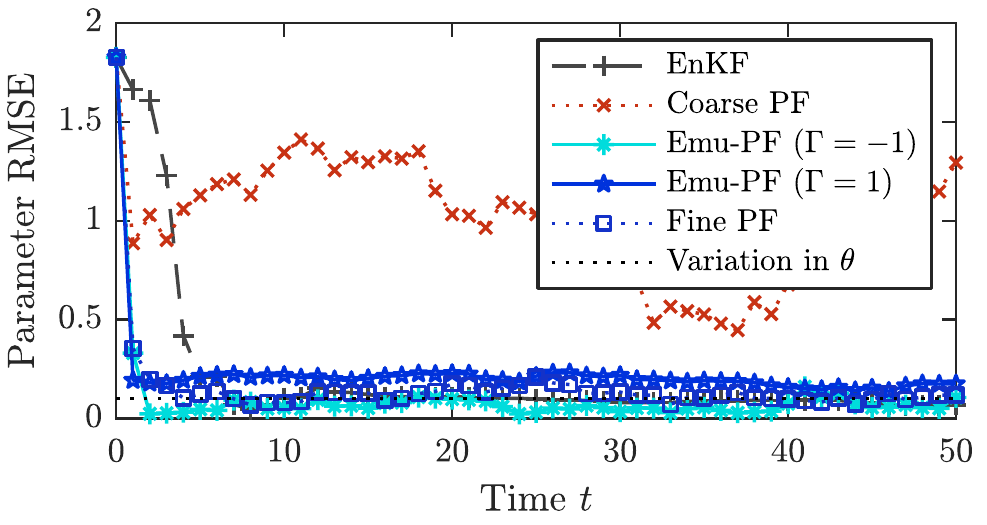}
     \includegraphics{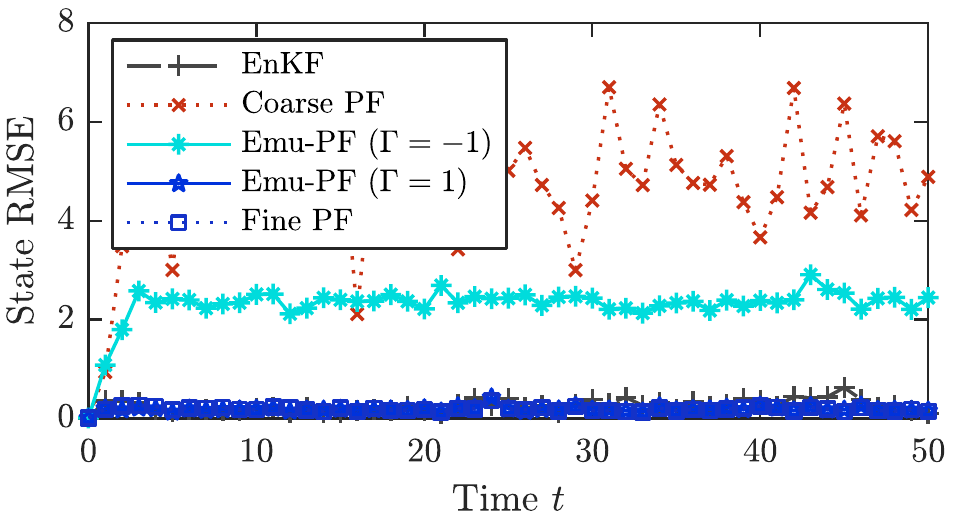}
     \caption{Error statistics for Experiment Three, $m=8$ observations at each observation time, of accuracy $\sigma_0=0.5$. The $\Gamma=-1$ Emu-PF provides accurate parameter estimates and reasonable state estimates, with a mean RMSE of $2.3$ for state variables (compared to $4.9$ for the coarse PF). However the $\Gamma=1$ Emu-PF is competitive with the much more expensive fine PF. }
     \label{fig:os5m8}
 \end{figure}

We now reduce the observation error $\sigma_o$ to $0.5$ and return to the fully observed $m=8$. In cases where the global dimension reduction of PCA is insufficient or impossible, perhaps because dimension reduction is only locally possible, then the localized inputs of \cref{sec:Epf} are an option. We test this approach with $\Gamma=1$, so input state variables are localized to three inputs, in \cref{fig:os5m8,table:os5m8}. 
 
 \begin{table}                                              
\centering                                                      
\begin{tabular}{ccccc}                                          
 & RMSE ($\theta$) & Var ($\theta$) & RMSE ($\bx$) & Var ($\bx$) \\
Fine PF & 0.11 & 0.0038 & 0.18 & 0.022 \\                       
Coarse PF & 0.84 & 0.0016 & 4.8 & 0.00032 \\                    
EnKF & 0.084 & 5.8e-06 & 0.3 & 0.0083 \\                        
Emu-PF ($\Gamma=-1$) & 0.062 & 0.00029 & 2.3 & 3 \\             
Emu-PF ($\Gamma=1$) & 0.18 & 0.0014 & 0.18 & 0.023           
\end{tabular}                                                   
\caption{Summary statistics for Experiment Three.}             
\label{table:os5m8}                                             
\end{table}    

\subsection*{Experiment Four: accurate, sparse observations} 
\begin{figure}
     \centering
     \includegraphics{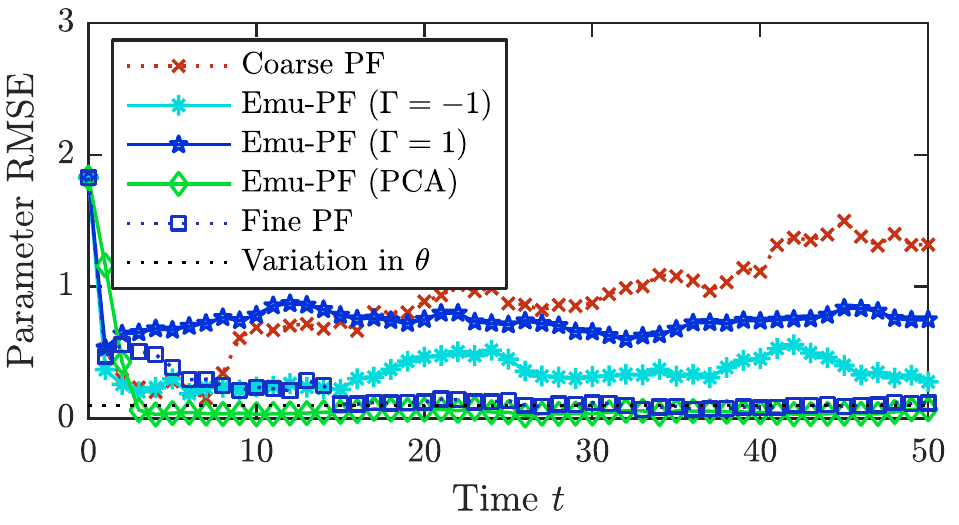}
     \includegraphics{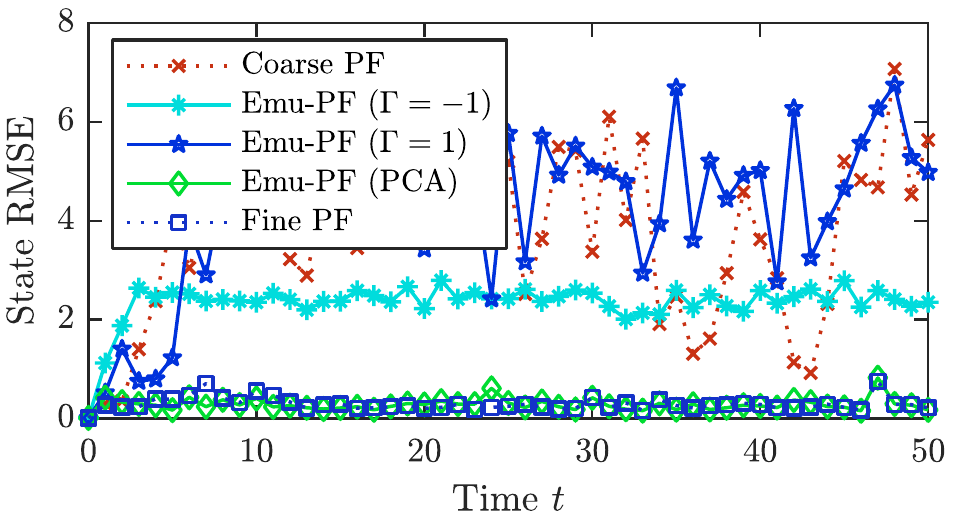}
     \caption{Error statistics for Experiment Four, $m=4$ observations at each observation time, of accuracy $\sigma_0=0.5$. In this case the $\Gamma=1$ Emu-PF performs only as well as the coarse PF. However the Emu-PF employing PCA is still competitive with the, much more expensive, fine PF.}
     \label{fig:os5m4}
 \end{figure}

We preserve $\sigma_o=0.5$ but reduce to $m=4$ observations. One drawback to the $\Gamma>0$ Emu-PF that we have observed is that it can be unstable if the filtering problem is slightly too hard; we infer that the emulator is given insufficient training data for the strongly localized input variables. \Cref{fig:os5m4,table:os5m4} show the $\Gamma=1$ Emu-PF performs significantly worse than the, technically inferior, Emu-PF with $\Gamma=-1$. In this case again the Emu-PF employing both PCA and PPE is competitive with the fine PF (though \cref{table:os5m4} shows the variance in the Emu-PF is smaller, indicating that it may be too tightly spread).
Localising strategies like the $\Gamma=1$ approach are critical in many modern DA applications. The results of Experiment Four demonstrate that the localization strategy we have adopted is insufficient for more difficult filtering problems. We plan for future work to combine such localization strategies with the dimension reduction strategy of \cref{sec:pcapf}.
 
\begin{table}                                                   
\centering                                                      
\begin{tabular}{ccccc}                                          
 & RMSE ($\theta$) & Var ($\theta$) & RMSE ($\bx$) & Var ($\bx$) \\
Fine PF & 0.095 & 0.0038 & 0.26 & 0.045 \\                      
Coarse PF & 1.1 & 0.0012 & 3.7 & 0.019 \\                       
EnKF & 0.11 & 2.4e-05 & 0.31 & 0.02 \\                          
Emu-PF ($\Gamma=-1$) & 0.38 & 0.00029 & 2.4 & 3 \\              
Emu-PF ($\Gamma=1$) & 0.72 & 0.00049 & 4.8 & 0.0034 \\          
Emu-PF (PCA) & 0.045 & 0.00096 & 0.26 & 0.034                 
\end{tabular}                                                   
\caption{Summary statistics for Experiment Four.}                
\label{table:os5m4}                                             
\end{table}

\subsection*{Experiment Five: state estimation}
\begin{figure}
    \centering
    \includegraphics{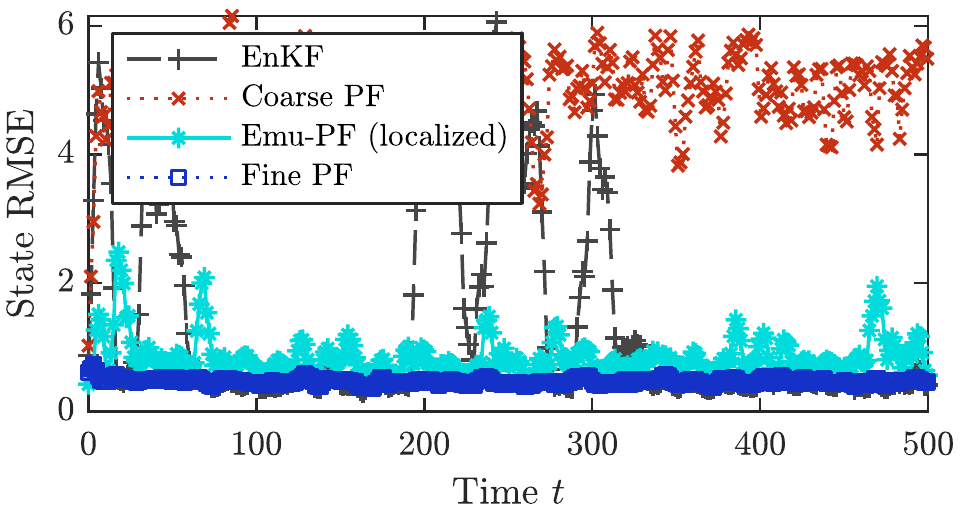}
    \caption{Summary statistics for Experiment Five, long-time state estimation with $m=4$ observations of accuracy $\sigma_o=1$. The median RMSE for EnKF and fine PF are similar; however the EnKF error occasionally spikes. The sliced Emu-PF of \cref{sec:Epf} is stable, with no large error spikes, and performs close to the fine PF in accuracy.}
    \label{fig:m4s}
\end{figure}

We showcase the localization strategy of \cref{sec:Epf}, the sliced Emu-PF with $\Gamma=1$. We fix $\theta=(2,\,1)^T$ in all methods, so that the only uncertainty is the state variables (but the filtering problem is still more difficult than the standard Lorenz-96, as the forcing \cref{forcing} is not uniform). In this state estimation experiment we assimilate every second variable with $m=4$ and standard accuracy $\sigma_0=1$, at each of $10,000$ observation times. Implementation for two algorithms differs in this experiment: the EnKF employs multiplicative inflation of $1.1$ (tuned to minimize RMSE) and the PF algorithms jitter particles with white noise of variance $0.01$ after each resampling step\footnote{this step is necessary here because the model is deterministic; if not jittered, the PF ensemble will collapse and all particles will be identical. Jittering is not strictly necessary for the PF in previous experiments because the parameter model is stochastic, and is not needed for the Emu-PF anywhere because the emulator already translates model uncertainty into noise.}. The Emu-PF with $\Gamma=1$ and the Emu-PF employing PCA (not plotted) both attained similar error values to the coarse PF\footnote{additionally, the Emu-PF with PCA halted due to an error with the PPE code.}. Results in \cref{fig:m4s,table:m4s} show the sliced Emu-PF outperform the EnKF and attain performance almost on par with the fine PF. 

\begin{table}                                   
\centering                                      
\begin{tabular}{ccc}                            
 & RMSE ($\bx$) & Var ($\bx$) \\                
Fine PF & 0.47 & 0.15 \\                        
Coarse PF & 5.1 & 0.16 \\                      
EnKF & 1 & 0.096 \\                             
Emu-PF (Localized) & 0.83 & 0.31 \\             
\end{tabular}                                   
\caption{Summary statistics for Experiment Five.}
\label{table:m4s}                               
\end{table}  



\section{Discussion and Future directions}
\label{sec:discussion}
In this work, we present a straight-forward utilization of statistical emulators within sequential data assimilation. We use random function models,  specifically Gaussian process emulators (GPs), to learn the mapping from state and/or parameter values at one observation instance to the next. This model-learning technique pairs well with particle filters that typically require $10^3-10^5$ forward model runs to assimilate each observation in time.
The gist of our methodology is that a GP provides interpolation between model forecasts -- thought of as functions of the parameter and/or previous state values at a fixed time -- and may be used to produce additional forecasts, and thus provide a cheap means to improve PF performance. Further, statistical emulators provide a built-in estimate of model performance in terms of the predictive variance of the Gaussian process. In our suite of simulation studies, we find that GP emulator-based particle filters with $100$ model runs perform on par or better when compared to a $10^4$ particle ``gold-standard" particle filter.

We explore several variations of the basic emu-PF algorithm, both to improve performance and to test various approaches to dimension reduction within the emulator. We introduce these various adaptions to mimic two salient flavors of dimension reduction on inputs to the dynamic forward mapping---namely two forms of localized dimension reduction, and a strategy for global dimension reduction. Localization is a widely-used and effective tool in DA to eliminate the impacts of long-range correlations on estimations and forecasts. The two approaches may be combined in future implementations of Emu-PF: one can imagine utilizing ``global" dimension reduction tools within the localization domain of a gridded model. We further utilize the parallel partial emulator in a variation of the emu-PF appropriate for functional or vector-valued model output.

We test this suite of algorithms through various simulation experiments on an 8-member Lorenz-96 system. We begin by considering a parameterized forcing that induces a bi-modal posterior distribution in parameter space. The emu-PF is able to obtain well-resolved bi-modal posteriors in parameter space with only $100$ forward-model runs. We then consider a series of assimilation experiments that present an increased challenge as we lower the dimension of the  observational space. We conclude that the success of the computationally cheap emu-PF with various forms of localization bodes well for this tool to be explored more widely. 

A very strong asset of this methodology is that it can readily be combined with other modern advances in sequential data assimilation. We intend to explore several of these combinations in future work. For example, the $40$-dimensional Lorenz-96 system could be simulated by employing the Optimal Proposal PF \cite{snyder2011particle} in conjunction with emulators (in addition to other modern PF innovations discussed in \cref{sec:algorithms}). Further,  the approach could be combined with the Localized PF \cite{Poterjoy16}. In this case, we envision a dimension reduction for the emulator based on the support of the localization(s) utilized within the Localized PF. The emulated-based particle filter also has the potential to work nicely with the Equal Weight PF \cite{van2010}. One can re-express the equivalent weights problem readily on the probability density functions obtained with emulators. Then one could sample from the resulting distribution. These advanced PF techniques devise approaches to overcome the challenge of searching large sample spaces; our contribution is effectively to accelerate the sampling procedure, so that more samples can be taken. 


\bibliographystyle{aims}
\bibliography{DA-master-2019}
\end{document}